\documentclass[12pt]{article}
\usepackage{epsfig}
\usepackage{graphicx}
\begin{document}

\title{\Large{\bf A Dynamic View of Circular Colorings}}

\author{Hong-Gwa Yeh\thanks{Partially supported by
National Science Council of R.O.C. under grant
NSC94-2115-M-008-015.}
\\
\normalsize   Department of Mathematics\\
\normalsize   National Central University\\
\normalsize   Jhongli City, Taoyuan 32001, Taiwan\\
\normalsize   National Center for Theoretical Sciences, Taiwan}

\date{\small  April 2006}

\maketitle

\newtheorem{theorem}{Theorem}
\newtheorem{lemma}[theorem]{Lemma}
\newtheorem{corollary}[theorem]{Corollary}
\newtheorem{fact}[theorem]{Fact}
\newtheorem{observation}[theorem]{Observation}
\newtheorem{example}[theorem]{Example}
\newcommand{\qed}{\hfill $\Box$ }
\newenvironment{proof}{
\par
\noindent {\bf Proof.}\rm}%
{\mbox{}\hfill\rule{0.5em}{0.809em}\par}
\def\l{\ell}
\def\V{V(\vec{G})}
\def\E{E(\vec{G})}
\def\G{\vec{G}}

\baselineskip=20pt
\parindent=1cm

\begin{abstract}
 The main contributions of this paper are three-fold. First, we use a
 dynamic approach based on Reiter's pioneering work
 on Karp-Miller computation graphs \cite{reiter} to give a
 new and short proof of Mohar's Minty-type Theorem \cite{mohar}.
 Second, we bridge circular colorings and discrete event dynamic
 systems to show that the Barbosa and Gafni's results
 on circular chromatic
 number \cite{bg,yeh-zhu} can be generalized to edge-weighted symmetric
 directed graphs.
 Third, we use the above-mentioned dynamic view of circular
 colorings to construct new improved lower bounds on the circular chromatic
 number of a graph. We show as an example that the circular
 chromatic number of the line graph of the Petersen graph can be
 determined very easily by using these bounds.
\end{abstract}


\section{Introduction}


 In this paper we discuss the possibilities to use dynamic
 techniques to investigate circular coloring problems.
 From these approaches we will see an intimate connection between
 circular colorings and discrete event dynamic systems.
 To explain our approaches, we first
 give some definitions from graph theory and scheduling that will
 be used in this paper. Throughout this paper all graphs are connected.
 An {\em edge-weighted directed graph} (or simply {\em edge-weighted digraph})
 is a directed graph $\vec{G}$
 together with an assignment $c$ of positive real weights to each
 directed edges and is denoted as $(\vec{G},c)$.
 We use $\V$ to denote
 the set of vertices and $\E$ to denote the set
 of directed edges. For simplicity of notation, the directed edge
 $(u,v)$ is written as $uv$ and is called an {\em arc}, the weight of the
 arc $uv$ in $(\vec{G},c)$ is written as $c_{uv}$. The arc $uv$ is called the
 {\em in-arc} of $v$ and {\em out-arc} of $u$.
 An edge-weighted digraph $(\vec{G},c)$ is said to be {\em symmetric}
 if arcs $uv$, $vu$ both exist
 and do not exist for all vertices $u,v$ in $\vec{G}$.
 A {\em directed cycle} (or simply {\em dicycle})
 $C$ of $(\vec{G},c)$ is a closed directed walk with
 no repeated vertices  except that the
 beginning vertex and ending vertex.

 To each arc $uv$ in $(\vec{G},c)$ we may assign a
 number $T_{uv}$ of {\em tokens} (a nonnegative integer).
 The function $T$ is called an {\em initial marking}
 of $(\vec{G},c)$ (or $\vec{G}$).
 An edge-weighted digraph $(\vec{G},c)$ equipped with an
 initial marking $T$ is denoted by $(\vec{G},c,T)$ and is called a
 {\em timed marked graph}.
 The {\em token count} (resp., {\em weight}) of a dicycle $C$
 in $(\vec{G},c,T)$ is defined as the value
 $\sum_{uv\in E(C)}T_{uv}$ (resp., $\sum_{uv\in E(C)}c_{uv}$)
 and is denoted by $|C|_T$ (resp., $|C|_c$).

 A time marked graph $(\vec{G},c,T)$ can be used to
 dynamically model the data
 movement in parallel computations.
 To model the dynamic behavior
 of data in the parallel computation environment, a vertex of a
 timed marked graph represents a task, an arc $uv$ represents a
 data channel with communication cost $c_{uv}$, and
 tokens on arc $uv$
 represent input data of the task vertex $v$.
 The edge weight $c_{uv}$ represents the time required by task
 vertex $u$ to place the result of its operation on arc $uv$
 (all the tokens on $uv$ are available as input data to the
 task node $v$). Dynamically, if task vertex $u$ operates at time $t$, then at time
 $t+ c_{uv}$, $u$ places one token on arc $uv$.

 A vertex $u$ is said to be {\em fireable} at time $t$ if each in-arc of $u$
 contains at least one token at time $t$.
 A vertex can start its operation
 only when it is fireable. When a vertex $u$ starts its operation we
 say that the vertex $u$ {\it fires}.
 Upon firing of vertex $u$, say at time $t$,
 $u$ removes one token from each of its in-arcs;
 if $uv$ is an out-arc of u, then at time $t+ c_{uv}$,
 $u$ places one token on arc $uv$. Let $T_{uv}(t)$
 denote the number of tokens on arc $uv$ at time $t$.
 With this notation, we have that $T(0)\equiv T$ is the initial marking
 of the time marked graph $(\vec{G},c,T)$
 (i.e., the distribution of tokens before the firings of
 vertices).
 During firings of vertices the distribution of
 tokens may change.
 A distribution of tokens in $\vec{G}$ is called
 a {\em marking}.
 Given a timed marked graph, in parallel computation we
 are interested in, for each vertex,
 determining a sequence of times called {\em schedule}
 at which the vertex
 starts its firing (operation) independently of what
 the other vertices are doing. These times must be such
 that each vertex, upon firing, is assured to be fireable.

 Formally, a {\em schedule} $f$ of $(\vec{G},c,T)$
 is a set of functions
 $\{f_u\}_{u\in \V}$, where function $f_u:\{1,2,3, \cdots\}\rightarrow
 {\cal R}$ having $f_u(k_1)<f_u(k_2)$ for $k_1<k_2$,
 where
 $f_u(k)=t$
 means that vertex $u$ starts its $k$th firing
 at time $t$ under the schedule $f$.
 A schedule $f$ of $(\vec{G},c,T)$ is called an {\em admissible schedule} if
 $(\vec{G},c,T)$ ``runs properly" according to the schedule $f$, i.e.,
 for every vertex $u$ and integer $k$,
 $f_u(k)=t$ implies $T^f_e(t)>0$ for all in-arc $e$ of $u$,
 where
 $T^f_e(t)$ is the number of tokens on arc $e$ at time $t$
 under the schedule $f$. Computer scientists are particularly
 interested in regular schedules. To that end, we make the
 following definition: An admissible schedule $f$ is
 {\em periodic with period} $p$ if there exist real numbers
 $x_u$ such that $$f_u(k)=x_u+p(k-1)$$
 for any vertex $u$ and positive integer $k$.

 Suppose $k\geq 2d$ are positive integers.
 A $(k,d)$-$coloring$ of an undirected graph  $G$ is a mapping $f:V(G)\rightarrow
 \{0,1,\ldots,k-1\}$ such that for any edge $xy$ of $G$,
 $d\leq |f(x)-f(y)|\leq k-d$. If $G$ has a $(k,d)$-coloring then we
 say $G$ is $(k,d)$-$colorable$.
 The {\em circular chromatic number} $\chi_c(G)$ of $G$ is defined
 as
 \begin{center}
 $\chi_c(G)=\inf\{k/d:$ $G$ is $(k,d)$-colorable $\}$.
 \end{center}
 In fact, to determine the circular chromatic number
 of $G$, it suffices to check finitely many $k$, $d$
 whether $G$ is $(k,d)$-colorable. In
 \cite{zhou,survey,survey2005} we see that
 $$\mbox{$\chi_c(G)\in \{ {k \over d} : \mbox{$k\leq |V(G)|$, $d\leq \alpha(G)$ and }
   { |V(G)| \over \alpha(G) } \leq {k \over d} \leq \chi(G) \}$,}$$
 where $\alpha(G)$ is the maximum size of an independent set in $G$.

 For a positive real $p$, let $S^p$ denote a circle with perimeter
 $p$ centered at the origin of ${\cal R}^2$. In the
 obvious way, we can identify the circle $S^p$ with the
 interval $[0, p)$. For $x,y\in S^p$,
 let $d(x,y)$ denote the length of the arc on $S^p$ from $x$ to $y$ in the
 clockwise direction if $x\not= y$, and let $d(x,y)=0$ if $x=y$.
 A {\em circular $p$-coloring} of an edge-weighted digraph $(\vec{G},c)$
 is a function $\varphi: \V\rightarrow S^p$ such that
 $d(\varphi(u),\varphi(v))\geq c_{uv}$ for each
 arc $uv$ in $\vec{G}$.
 The
 {\em circular chromatic number $\chi_c(\vec{G},c)$ of an
 edge-weighted digraph} $(\vec{G},c)$,
 recently introduced by Mohar \cite{mohar}, is defined as
 \begin{center}
 $\chi_c(\vec{G},c)=\inf\{p:$ $(\vec{G},c)$ has a circular
 $p$-coloring $\}$.
 \end{center}
 It was shown in \cite{mohar} that
 the notion of $\chi_c(\vec{G},c)$
 generalizes several well-known optimization problems,
 such as the
 usual circular chromatic number \cite{survey},
 the weighted circular colorings \cite{deuber},
 the linear arboricity of a graph and
 the metric traveling salesman problem.


\section{Circular coloring and
periodic admissible schedule}


 In this section, we will use a dynamic view of the
 circular chromatic number $\chi_c(\vec{G},c)$ (which is defined
 below) to give a
 new and short proof of Mohar's Minty-type Theorem \cite{mohar}.
 The key techniques we use here
 are Theorems \ref{cc-pas} and \ref{reiter}.
 Theorem \ref{cc-pas} shows
 a relationship between circular colorings and
 periodic schedules of timed marked graphs.
 Theorem \ref{reiter} is the Reiter's theorem on the
 period of a periodic admissible schedule.
 An initial marking $T$ of an edge-weighted symmetric digraph
 $(\vec{G},c)$ is said to be {\em good} if
 for each dicycle $C$ and each arc $uv$ in $\vec{G}$,
 $$ |C|_T>0 \mbox{\,\,\,and\,\,\,} T_{uv}+T_{vu}=1.$$
 We remark that
 a good initial marking $T$ of $(\vec{G},c)$ is said to be {\it live}
 and {\it safe} in the terminology of \cite{marked-graphs}.
 A marking $T$ of $\vec{G}$ is called {\em live} (see \cite{marked-graphs})
 if every vertex of
 $\vec{G}$ is fireable upon this token distribution $T$, or can be
 made fireable through some sequence of firings. A marking $T$ is
 said to be
 {\em safe} (see \cite{marked-graphs}) if the token count of an arc never exceeds 1.



\begin{theorem}
\label{cc-pas}
 An edge-weighted symmetric digraph $(\vec{G},c)$
 has a circular $p$-coloring if and only if
 there is a good initial marking $T$ for $(\vec{G},c)$
 for which the timed marked graph $(\vec{G},c,T)$
 admits a periodic admissible schedule with period $p$.
\end{theorem}

\begin{proof}
 $(\Rightarrow)$ Let $\varphi$ be a circular $p$-coloring
 of $(\vec{G},c)$. Define an initial marking $T$ for
 $(\vec{G},c)$ as follows: for each arc $uv$ of $\vec{G}$,
 $$\left\{
 \begin{array}{l}
    \mbox{$T_{uv}=0$ and $T_{vu}=1$,\,\,\,\, if $\varphi(v)> \varphi(u)$;}\\
    \mbox{$T_{uv}=1$ and $T_{vu}=0$,\,\,\,\, otherwise.}
 \end{array}\right.
 $$
 We claim that $T$ is a good initial marking $T$
 for the edge-weighted digraph $(\vec{G},c)$.
 Indeed, by the choice of $T$, any
 zero token count dicycle $C=u_1 u_2 \cdots u_k u_1$
 would lead to a contradiction $\varphi(u_1)<\varphi(u_1)$.
 Next, let us construct a schedule $f$
 for the timed marked graph $(\vec{G},c,T)$
 as follows: for each vertex $u$ and positive integer $k$,
 define
 $$ f_u(k)=\varphi(u)+p(k-1).
 $$
 The fact that
 $d(\varphi(u),\varphi(v))\geq c_{uv}$ for each
 arc $uv$ of $\vec{G}$ immediately implies that
 $f=\{f_u\}_{u\in \V}$ is an admissible schedule for $(\vec{G},c,T)$.
 Here we note that $T$ is a live and safe marking.

 $(\Leftarrow)$ Suppose that $(\vec{G},c,T)$ has a periodic
 admissible schedule $\varphi$ having $ \varphi_u(k)=x_u+p(k-1)  $
 for each vertex $u$ and positive integer $k$, where $x_u\geq 0$.
 For any arc $uv$, say $T_{uv}=1$ and $T_{vu}=0$,
 since $\varphi$ is an admissible schedule for $(\vec{G},c,T)$,
 it must be
 $$\left\{
 \begin{array}{lcl}
    \varphi_u(1)& \geq & \varphi_v(1)+c_{vu};\\
    \varphi_v(2)& \geq & \varphi_u(1)+c_{uv},
 \end{array}\right.
\mbox{\,and hence\,\,}
 \left\{
 \begin{array}{lcl}
    x_u   & \geq & x_v+c_{vu};\\
    x_v+p & \geq & x_u+c_{uv}.
 \end{array}\right.
 $$
 That is $c_{vu}\leq x_u-x_v\leq p-c_{uv}$.
 It turns out that
 there is an obvious way to define
 a circular $p$-coloring $\varphi$
 for $(\vec{G},c)$.
 \end{proof}

 Theorem \ref{reiter} below, due to Reiter's pioneering work
 \cite{reiter} on Karp-Miller computation graphs \cite{km}, gives a
 necessary and sufficient condition for the existence of a periodic
 admissible schedule with specific period.
 We note that similar results had been presented in
 \cite{murata,ho,timed-petri-nets}.
 We say that vertex $u$ can {\em reach}
 vertex $v$ in digraph $\vec{G}$
 if there is a directed path from $u$ to $v$.



\begin{theorem}
\label{reiter}
 Suppose $(\vec{G},c,T)$ is a timed marked graph with the
 following two properties: There is a vertex $s$ in $\vec{G}$ which can reach
 every other vertex; Each dicycle $C$ of $\vec{G}$ has $|C|_T>0$.
 Then $(\vec{G},c,T)$
 has a periodic admissible schedule with period $p$
 if and only if
 $$
 p\geq \max_{C}{|C|_c\over |C|_T},
 $$
 where the maximum is taken over all dicycles $C$ in $\vec{G}$.
\end{theorem}


 By using Hoffman's circulation theorem,
 in \cite{mohar} Mohar proved Theorem \ref{mohar} stated below.
 This result is a
 Minty-type theorem for
 circular chromatic number of an edge-weighted digraph and
 generalizes the corresponding results in
 \cite{gtz,minty}.



 \begin{theorem}
 \label{mohar}
 Let $(\vec{G},c)$ be an edge-weighted symmetric digraph.
 Then
 $$
 \chi_c(\vec{G},c)=\min_{T}\max_{C}{|C|_c\over |C|_T},
 $$
 where
 the minimum is taken over all good initial markings $T$ of $(\vec{G},c)$
 and
 the maximum is taken over all dicycles $C$ of $\vec{G}$.
 \end{theorem}

 It is easy to see that Mohar's Minty-type theorem follows from
 Theorems \ref{cc-pas} and  \ref{reiter}  immediately.
 Our approach here provides a new insight
 into the circular chromatic numbers.
 It connects several different research fields,
 timed Petri nets \cite{timed-petri-nets},
 marked graphs \cite{marked-graphs},
 scheduling,
 and coloring.
 It reveals that there are strong coherent and coordinated
 relationships between these topics.
 We feel sure that the results in this paper can shed insight on
 these areas.


 In the rest part of this section, we turn our attention to
 edge-weighted symmetric digraphs $(\vec{G},c)$ with weights
 $c:E(\vec{G})\rightarrow R^+\cup \{0\}$
 such that $c_{uv}+c_{vu}\not =0$ for each arc $uv$ in $\vec{G}$.
 A {\em weak circular $p$-coloring} of an edge-weighted
 (not necessarily symmetric) digraph $(\vec{G},c)$
 is defined to be a function $\varphi: \V\rightarrow S^p$
 which satisfies the following two conditions:
 (1) For each arc $uv$ in $\vec{G}$,
 either $\varphi(u)=\varphi(v)$ or $d(\varphi(u),\varphi(v))\geq
 c_{uv}$.
 (2) For any point $x$ of $S^p$, the subdigraph of $\vec{G}$ induced by
 $\varphi^{-1}(x)$ contains no dicycles $C$ having $c_{uv}>0$ for
 each arc $uv$ in $\vec{G}$.

 The weak circular coloring of an edge-weighted digraph
 defined here generalizes the earlier definition of weak circular
 coloring of digraphs in \cite{bfjkb}, and
 it is easy to see that
 $$\chi_c(\vec{G},c)=\inf\{p\geq\max_{uv\in E(\vec{G})}(c_{uv}+c_{vu}):
 (\vec{G},c)\mbox{ has a weak circular $p$-coloring}\}.$$
 A weak circular coloring $\varphi$ of
 an edge-weighted symmetric digraph $(\vec{G},c)$ is said to have
 a {\em tight dicycle} $C$ if for each arc $uv$ of $C$,
 $d(\varphi(u),\varphi(v))=c_{uv}$.
 The following theorem strengthens and generalizes Theorem 3.3 of
 \cite{survey2005}, and thus is both a generalization and a strengthening of
 the corresponding results in \cite{guichard, mohar}.



 \begin{theorem}
 \label{tight}
 Let $(\vec{G},c)$ be an edge-weighted symmetric digraph with
 $c:E(\vec{G})\rightarrow R^+\cup \{0\}$
 such that $c_{uv}+c_{vu}\not =0$ for each arc $uv$ in $\vec{G}$.
 Suppose $r$ is greater than but not equal to $\max_{uv\in E(\vec{G})}(c_{uv}+c_{vu})$.
 Then $\chi_c(\G,c)=r$ if and only if there is a weak circular
 $r$-coloring of $(\G,c)$, and  every weak circular
 $r$-coloring of $(\G,c)$ has a tight dicycle.
 \end{theorem}
 \begin{proof}
 $(\Rightarrow)$
  Suppose that $\chi_c(\G,c)=r$.
  By Theorem \ref{mohar} there is a good initial marking $T$ of $(\vec{G},c)$
  such that $$r=\max_{C}{|C|_c\over |C|_T},$$
 where
 the maximum is taken over all dicycles $C$ of $\vec{G}$.
 Define a mapping $w:E(\G)\rightarrow R$ such that
 $w(uv)=c_{uv}-rT_{uv}$ for each arc $uv$ in $(\G,c)$.
 We'll abbreviate $w(uv)$ to $w_{uv}$ from now on.
 Next, for each vertex $u$ in $\G$,
 let $\varphi(u)$ denote the maximum weight
 (with respect to the arc weights $w_{uv}$) of
 a directed walk of $\G$ ending at $u$.
 Our aim is to show that a weak circular
 $r$-coloring of $(\G,c)$ can be obtained from $\varphi$.

 Since for each dicycle $C$ of $\G$
 $$ |C|_w=|C|_c-r|C|_T=|C|_T\left({|C|_c\over |C|_T}-r\right)\leq 0,
 $$
 $\varphi(u)$ is well-defined for each vertex $u$.
 From the definition of $\varphi$ it is easy to see that
 if $uv$ is an arc of $\G$, then we have
 $$w_{uv}\leq \varphi(v)-\varphi(u)\leq -w_{vu}.
 $$
 Denote by $V_u$ the vertex subset $\{x\in V(\G): \varphi(x)\equiv
 \varphi(u)\mbox{ (mod $r$)}\}$. We make the following claim.

 \noindent
 {\bf Claim} {\em For any vertex $u$ in $\G$, the subgraph of $(\G,c)$ induced by $V_u$ contains no
 dicycles $C$ having $c_{xy}>0$ for each arc $xy$ in $C$.}

 To prove this claim by contradiction, suppose that the claim is
 false, and let $C$ be a dicycle of $(\G,c)$ such that
 $\varphi(x)=\varphi(y)\mbox{ (mod $r$)}$
 for any two vertices $x,y$ in $C$,
 and $c_{uv}>0$ for each arc $uv$ in $C$.
 Without loss of generality, let $C=(1,2,3,\ldots,t,1)$ and
 $\varphi(1)=\min_{1\leq i \leq t}\varphi(i)$.
 If $\varphi(t)>\varphi(1)$ then for some positive integer $k$ we
 have $\varphi(t)=\varphi(1)+kr$, and hence
 $$ \varphi(1)\geq \varphi(t)+w_{t1}=\varphi(1)+kr+c_{t1}-rT_{t1}.
 $$
 It follows that $0\geq c_{t1}+(k-T_{t1})r\geq c_{t1}>0$, a
 contradiction.
 Hence it must be $\varphi(t)=\varphi(1)$.
 Similar arguments yield that
 $\varphi(1)=\varphi(2)=\varphi(3)=\cdots=\varphi(t)$.
 Because $\varphi(2)\geq \varphi(1)+w_{12}$, we see that
 $0\geq c_{12}-rT_{12}$ And so $T_{12}=1$, since $c_{12}>0$ and
 $T_{12}\in\{0,1\}$.
 By the same argument, we have
 $$T_{12}=T_{23}=T_{34}=\cdots
 =T_{(t-1)t}=T_{t1}=1,$$
 and thus $|\hat{C}|_T=0$ where $\hat{C}$ is the dicycle $(t,t-1,t-2,\ldots,2,1,t)$ of $\vec{G}$,
 that is impossible, since $T$ is a good initial marking of $(\vec{G},c)$.
 This completes the proof of the claim.

 We then define a mapping $\hat{\varphi}$ from $V(\vec{G})$ to
 $S^{r}$
 such that $\hat{\varphi}(v)=\varphi(v)({\rm mod}\, r)$ for each vertex $v$ of $V(\vec{G})$.
 It is easy to verify that $\hat{\varphi}$ is a weak circular $r$-coloring of
 $(\vec{G},c)$.

 Next, let $g$ be an arbitrary weak circular $r$-coloring of
 $(\vec{G},c)$. We are going to show that $g$ has a tight dicycle.
 Assume that $g$ has no tight dicycles.
 We construct a digraph $D(\vec{G},c,g)$ as follows: It has vertex
 set $V(\vec{G})$, and there is an arc from $x$ to $y$ in $D(\vec{G},c,g)$
 if and only if $xy$ is an arc of $\vec{G}$
 that satisfies one of the two conditions below,

 (1) $g(x)=g(y)$ and $c_{xy}=0$;

 (2) $g(x)\not =g(y)$ and $d(g(x),g(y))=c_{xy}$.

  Since $g$ has no tight dicycles, $D(\vec{G},c,g)$ contains no
  dicycles.
  Suppose $s$ is a sink (that is a vertex of outdegree 0) with indegree at least $1$
  in the acyclic digraph $D(\vec{G},c,g)$.
  Clearly, there exists sufficiently small $\epsilon >0$  such that
  the mapping $g':V(\vec{G})\rightarrow S^r$ defined below is a weak circular $r$-coloring of
  $(\vec{G},c)$, and $E(D(\vec{G},c,g'))$ is a proper subset of $E(D(\vec{G},c,g))$.
  $$ g'(v)=\left\{\begin{array}{lll}
                 g(v)+\epsilon & \mbox{if} & v=s;\\
                 g(v)          & \mbox{if} & v\not= s.
                 \end{array}\right. $$
  We may repeat this process many times to obtain a weak circular $r$-coloring $\hat{g}$ of
  $(\vec{G},c)$ such that $D(\vec{G},c,\hat{g})$
  contains no arcs.
  Let $g'$ be the mapping from $V(\vec{G})$ to $S^{r-\delta}$
  defined by $g'(v)=(1-\delta/r)\hat{g}(v)$ for any $v\in V(\vec{G})$.
  It is easy to see that for sufficiently small $\delta >0$,
  $g'$ is a weak circular $(r-\delta)$-coloring of $(\vec{G},c)$,
  and we arrive at a contradiction. Hence
  $g$ must have a tight dicycle.

  $(\Leftarrow)$ To prove this direction it is enough to show that
  $\chi_c(\vec{G},c)\geq r$. Assume $\chi_c(\vec{G},c)=r' <r$.
  From the first part of the proof, we see that there is a
  weak circular $r'$-coloring $\varphi '$ of $(\vec{G},c)$.
  Let $\varphi$ be a mapping from $V(\vec{G})$ to $S^r$ such that
  $\varphi(v)={r\over r'}\varphi'(v)$.
  It is clear that $\varphi$ is a weak circular $r$-coloring of
  $(\vec{G},c)$, and $\varphi$ has no tight dicycles.
  That is a contradiction, and hence $\chi_c(\vec{G},c)=r$.
  \end{proof}


 \section{Circular coloring and
 discrete event dynamic system}

 In this section
 all edge weights are assumed to be positive integers.
 It is not hard to see that
 all results in this section are also
 true for positive rational edge weights.

 Suppose we are playing
 a token game in
 a connected edge-weighted symmetric digraph
 $(\vec{G},c)$ which
 has a good initial marking
 $T$.
 We assume that
 the time pulse $T_{pulse}$ is discrete, say
 $T_{pulse}\in \{1,2,3,\cdots\}$.
 Starting from $T_{pulse}=1$, at each time pulse
 we fire all fireable vertices simultaneously.
 A set of functions $g=\{g_u\}_{u\in \V}$ is defined
 to record the changes of token distributions in $\vec{G}$
 such that
 $g_u(k)=t$ means that vertex $u$ starts its $k$th firing
 at $T_{pulse}=t$ during the game.
 This set of functions $g$ can serve as
 an admissible schedule of $(\vec{G},c,T)$.
 Define a mapping $h:\{1,2,3,\cdots \}\rightarrow 2^{\V}$
 such that
 $h(t)$ consists of
 the vertices that are fireable at time pulse $T_{pulse}=t$.
 Here $h(t)$ may be an empty set due to the positive
 edge weights (time delays).
 As time goes on,
 the pigeon-hole principle and
 the fact that all edge weights are integral
 ensure that the dynamic system will converge to
 the following steady state:
 There exist
 positive integers $p$ and $M$
 such that
 $$ h(t+p)=h(t)
 $$
 as $t\geq M$.
 The value of such $p$ at its minimum is denoted by
 $p_{(\vec{G},c,T)}$ and is called the {\em period} of $(\vec{G},c,T)$.
 Let $m^u_{(\vec{G},c,T)}$
 denote the number of times
 that vertex $u$ becomes fireable in the time interval $[M,M+p)$.
 By \cite{bg-proceeding,bg,interleaved}, it can be shown that
 $$m^u_{(\vec{G},c,T)}=m^v_{(\vec{G},c,T)}
 $$
 for any two vertices $u$ and $v$,
 and hence we can omit the superindex $u$
 from $m^u_{(\vec{G},c,T)}$.
 The integer $m_{(\vec{G},c,T)}$ is called
 the {\em multiplicity} of $(\vec{G},c,T)$.
 With these notations we prove
 the following characterization of steady state
 of the above token game on $(\vec{G},c,T)$.



 \begin{lemma}
 \label{steady-state}
 $${p_{(\vec{G},c,T)}\over m_{(\vec{G},c,T)}}
  =
 \max_{C}{|C|_c\over |C|_T},
 $$
 where the maximum is taken over all dicycles $C$ of $\vec{G}$.
 \end{lemma}

 \begin{proof}
 Let $p= \max_{C}{|C|_c\over |C|_T}$.
 Theorem \ref{reiter} says that
 $(\vec{G},c,T)$ has a periodic admissible schedule $f$
 of the form $f_u(k)=x_u+p(k-1)$,
 where $x_u$'s are non-negative real numbers.
 Let $\widetilde{f}_u(t)$ (resp. $\widetilde{g}_u(t)$)
 denote the number of fires of vertex $u$
 during the time interval $[0,t]$
 under the schedule $f$ (resp. $g$).
 The following intuitive result
 can be proved by induction on the time pulse $t$.\\

 {\bf Claim:}
 $\widetilde{g}_u(t) \geq \widetilde{f}_u(t)$
 {\em for any vertex $u$ and any time pulse $t$.}\\

 \noindent
 Let $\widehat{C}$ be a dicycle of $\vec{G}$
 with $|\widehat{C}|_c/|\widehat{C}|_T=p$.
 Let $u$ be a vertex on $\widehat{C}$.
 Since the token distribution changes over time during playing the
 token game on $(\vec{G},c,T)$
 and eventually converges to a steady state with
 period $p_{(\vec{G},c,T)}$ and
 multiplicity $m_{(\vec{G},c,T)}$, we have
 \begin{eqnarray*}
 {  p_{(\vec{G},c,T)}  \over  m_{(\vec{G},c,T)}  }
 &   =  &
 \lim_{t\rightarrow \infty} {  t \over  \widetilde{g}_u(t) }\\
 & \leq &
 \limsup_{t\rightarrow \infty} {  t \over  \widetilde{f}_u(t) }\\
 &   =  &
 \limsup_{k\rightarrow \infty} {  f_u(k) \over  k             }\\
 &   =  &
 \limsup_{k\rightarrow \infty} {  x_u+p(k-1) \over  k         }\\
 &   =  & p.
 \end{eqnarray*}
On the other hand,
 \begin{eqnarray*}
 {  p_{(\vec{G},c,T)}  \over  m_{(\vec{G},c,T)}  }
 &   =  &
 \lim_{t\rightarrow \infty} {  t \over  \widetilde{g}_u(t) }\\
 &   \geq  &
 \liminf_{t\rightarrow \infty} {  t|\widehat{C}|_c \over  \widetilde{g}_u(t|\widehat{C}|_c) }\\
 & \geq &
 \liminf_{k\rightarrow \infty} {  t|\widehat{C}|_c \over  (t+1)|\widehat{C}|_T    }\\
 &   =  & p.
 \end{eqnarray*}
 \end{proof}

 By Theorem \ref{mohar} and Lemma \ref{steady-state} we have
 the following ``dynamic" characterization of circular chromatic
 number.



 \begin{theorem}
 \label{dynamic}
 Let $(\vec{G},c)$ be an integral edge-weighted symmetric digraph.
 Then
 $$
 \chi_c(\vec{G},c)=\min_{T}{p_{(\vec{G},c,T)}\over m_{(\vec{G},c,T)}},
 $$
 where
 the minimum is taken over
 all good initial markings $T$ of $(\vec{G},c)$.
 \end{theorem}

 Let $\omega$ be an acyclic orientation
 of a connected undirected simple graph $G$.
 We denote by $(G,\omega)$ the corresponding oriented graph.
 Sometimes, by abuse of notation,
 we will simply write $\omega$ instead of $(G,\omega)$.
 A vertex in $(G,\omega)$ with zero outdegree (resp., indegree)
 is called a {\em sink} (resp., {\em source}).
 Let sink$(\omega)$ (resp., source$(\omega)$)
 denote the set of sinks (resp., sources) in $(G,\omega)$.
 One can obtain a sequence of acyclic orientations
 $\omega_0$, $\omega_1$, $\omega_2,\cdots$
 from $\omega$
 in such a way that  $\omega_0=\omega$,
 and $\omega_i$ ($i\geq 1$) is obtained from $\omega_{i-1}$
 by reversing the orientations of the edges incident to the sinks
 of $w_{i-1}$. We say that $\{\omega_i\}_{i=0}^\infty$ is {\em generated by} $\omega$.
 In Figure 1 we depict a sequence of acyclic
 orientations $\{\omega_i\}_{i=0}^\infty$ on $C_5$
 which is generated by $\omega_0$. Note that this sequence has
 the periodic property: $\omega_i=\omega_{i+5}$ for any $i\geq 0$.



 \begin{center}
 \qquad\includegraphics[scale=0.7]{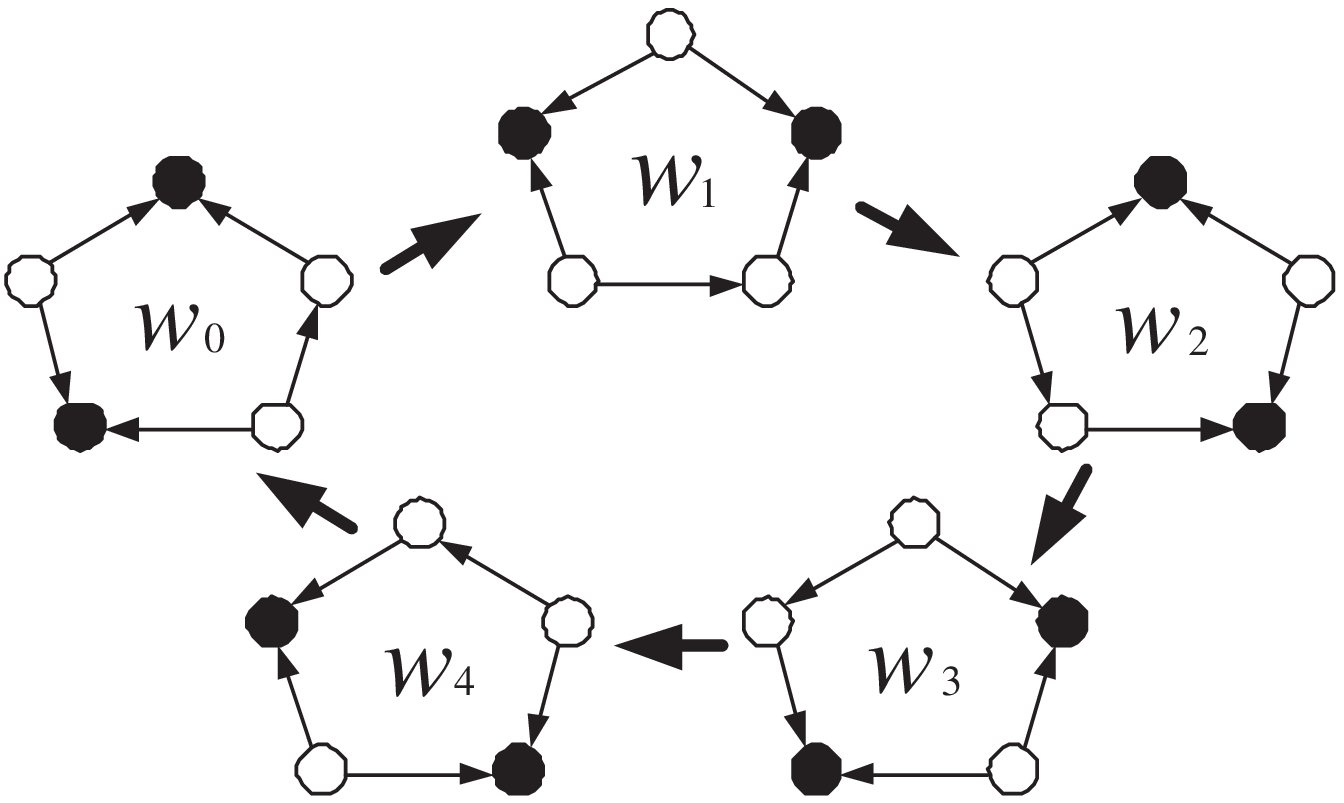}\\~\\
  Figure 1:
  A sequence of acyclic
  orientations $\{\omega_i\}_{i=0}^\infty$
  generated by $\omega_0$
 \end{center}

 Obviously the sequence $\{\omega_i\}_{i=0}^\infty$ generated by an acyclic orientation $\omega$
 has the following periodic behavior:
 There exist positive integers $p$ and $M$ such that
 $\omega_i=\omega_{i+p}$ for every $i\geq M$.
 The value $p$ at its minimum is denoted by $p_\omega$
 and is called the {\em period} of $\omega$.
 The $p_\omega$ acyclic orientations
 $\omega_M$, $\omega_{M+1},\cdots,\omega_{M+p_\omega-1}$
 is called a {\em period generated by} $\omega$.
 For any $k \geq M$ and $\ell \leq p_\omega -1$, the sequence
 $$\mbox{$\langle
 |$sink$(\omega_k)|$,$|$sink$(\omega_{k+1})|$,$\cdots$,$|$sink$(\omega_{k+\ell})|\rangle$}$$
 is called a {\em sub-pattern
 of a period generated by} $\omega$,
 sometimes we call it a {\em sub-pattern of $\omega$} for short.
 Let $m_\omega^u$ denote the number of times that the
 vertex $u$ becomes a sink in the period generated by $\omega$.
 It was shown in \cite{bg} that
 $m_\omega^u=m_\omega^v$ for any two vertices $u$ and $v$.
 So we write $m_\omega$ instead of $m_\omega^u$, and
 $m_\omega$ is called the {\em multiplicity} of $\omega$.
 Clearly, the following result, shown in \cite{bg,yeh-zhu},
 is a special case of Theorem \ref{dynamic}.



 \begin{corollary}
 \label{pm}
 Let $G$ be a connected simple graph.
 Then
 $$
 \chi_c(G)=\min_{\omega}{p_\omega \over m_\omega},
 $$
 where
 the minimum is taken over
 all acyclic orientations $\omega$ of $G$.
 \end{corollary}


 \section{Lower bounds for $\chi_c$}

 In this section, new lower bounds on the circular chromatic number
 of an undirected graph are derived by using
 the above-mentioned dynamic view of $\chi_c$.
 With the aid of these lower bounds,
 we show that circular chromatic numbers
 of some specific graphs can be
 determined by using a unified methodology.
 Note that
 these chromatic numbers were determined by using different ad hoc
 methods.
 In this section, for an acyclic orientation $\omega$ of a graph $G$,
 we are going to study the pattern of the sequence
 $$ \mbox{$|\rm{sink}(\omega_1)|$,$|\rm{sink}(\omega_2)|$,
 $|\rm{sink}(\omega_3)|$,\ldots,$|\rm{sink}(\omega_{p_\omega})|$, } $$
 where $\omega_1,\omega_2,\cdots,\omega_{p_\omega}$
 is the period generated by $\omega$.
 The study aims to provide interesting lower bounds for $\chi_c(G)$.
 The results in this section will show that
 the dynamic approaches in the previous sections
 do provide interesting insight into
 circular chromatic numbers.

 To simplify our expressions, throughout this section we assume
 that if the sequence of acyclic orientations
 $\omega_1,\omega_2,\cdots,\omega_{p_\omega}$
 is a period generated by $\omega$, then we automatically have
 $\omega_i=\omega_{i+p_\omega}$ for any integer $i\geq 1$.
 For a vertex $u$ of a graph $G$,
 let $N_k(u)$ denote all vertices of distance
 $k$ from $u$ in $G$, i.e.  $N_k(u)=\{v\in V(G): d_G(u,v)=k\}$.
 For a set $S\subseteq V(G)$, we define
 $N_1(S)=\{u\in V(G)\setminus S: uv\in E(G)$ for some $v\in S\}$.
 We write $N_1(u)$ (resp. $N_1(u,v)$) instead of
 $N_1(\{u\})$ (resp. $N_1(\{u,v\})$) for short.
 We say that a graph $G$ is $k$-{\it colorable} if $V(G)$ can be
 colored by at most $k$ colors so that adjacent vertices are colored
 by different colors. The {\it chromatic number} of $G$, denoted by
 $\chi(G)$, is the smallest $k$ such that $G$ is $k$-colorable.
 Let $\alpha_k(G)$ (or simply $\alpha_k$ if it cause no confusion)
 denote the maximum number of vertices
 in a vertex-induced $k$-colorable subgraph of $G$.
 Notice that $\alpha_1(G)=\alpha(G)$.
 For a set $S\subseteq V(G)$, by abuse of notation, we also use $S$ for the
 subgraph
 of $G$ induced by $S$.

 The following two theorems reveal the connection between the circular
 chromatic number of a graph and the chromatic number of its
 vertex's
 distance-1 neighborhood.



 \begin{theorem}
 \label{d1}
 For any vertex $u$ of a graph $G$,
 $\chi_c(G)\geq \chi(N_1(u))+1$.
 \end{theorem}

 \begin{proof}
 Let $\xi=\chi(N_1(u))$. By Corollary \ref{pm}, there is an acyclic
 orientation $\omega$ of $G$ such that $p_\omega/m_\omega=\chi_c(G)$.
 Let $\omega_1,\omega_2,\cdots,\omega_{p_\omega}$ be the period
 generated by $\omega$.
 For a vertex $u$, let $I_i^u$ to be $1$ if $u\in {\rm sink}(\omega_i)$ and
 $0$ otherwise. Note that
 $$ \sum_{i=1}^{p_\omega} I_i^u=m_\omega.$$
 If $u\in {\rm sink}(\omega_i)$ then clearly $$u\not\in
 \bigcup_{s=1}^\xi {\rm sink}(\omega_{i+s}).$$
 It follows that $p_\omega\geq
 \sum_{i=1}^{p_\omega}(\xi+1)I_i^u=(\xi+1)m_\omega$,
 and hence $p_\omega/m_\omega\geq \xi+1$.
 \end{proof}



 \begin{theorem}
 \label{alphat} Let $t$ be a positive integer.
 If every vertex $u$ of $G$ has $\chi(N_1(u))\geq t-1$,
 then $\chi_c(G)\geq {t|V(G)| /\alpha_t(G)}$.
 \end{theorem}

 \begin{proof}
 Let $\omega$ be an acyclic
 orientation of $G$ such that $p_\omega/m_\omega=\chi_c(G)$.
 Let $\omega_1,\omega_2,\cdots,\omega_{p_\omega}$ be the period
 generated by $\omega$.
 Since every vertex $u$ of $G$ has $\chi(N_1(u))\geq t-1$,
 we have $$ \mbox{${\rm sink}(\omega_{i})\bigcap{\rm
 sink}(\omega_{i+\ell})=\emptyset$,} $$
 for any index $i$ and any integer $\ell\in \{ 1,2,\cdots,t-1\}$.
 This implies that the inequality below holds:
 $$ t|V(G)|m_\omega=\sum_{i=1}^{p_\omega}\sum_{j=0}^{t-1}|{\rm sink}(\omega_{i+j})|
 \leq \sum_{i=1}^{p_\omega} \alpha_t(G)=p_\omega \alpha_t(G).
 $$
 Thus $p_\omega/m_\omega\geq {t|V(G)| /\alpha_t(G)}$.
 \end{proof}

 By Theorem \ref{d1}, one can easily show that the odd wheel depicted in
 Figure 2 has circular chromatic number 4.


 \begin{center}
 \qquad\includegraphics[scale=0.7]{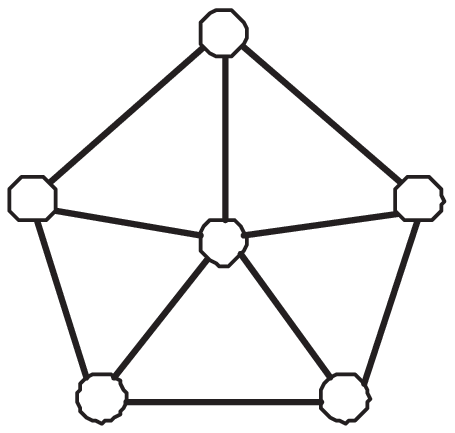}\\~\\
  Figure 2: An odd wheel.
 \end{center}

 Note that that Theorem \ref{alphat} says that every connected graph $G$
 has the property $\chi_c(G)\geq {2|V(G)| /\alpha_2(G)}$. This lower
 bound is sharp in the sense that there are infinitely many graphs
 $G_n$ that attain this bound. The graph $G_n$ is obtained from the
 even cycle $C_{8n}=(v_0,v_1,\ldots,v_{8n-1}, v_0)$ by joining $v_{2i}$
 and $v_{2i+4n}$ for $i=0,1,2,\ldots,2n-1$.
 As an example, the graph $G_{2}$ is depicted in Figure 3.


 \begin{center}
 \qquad\includegraphics[scale=0.4]{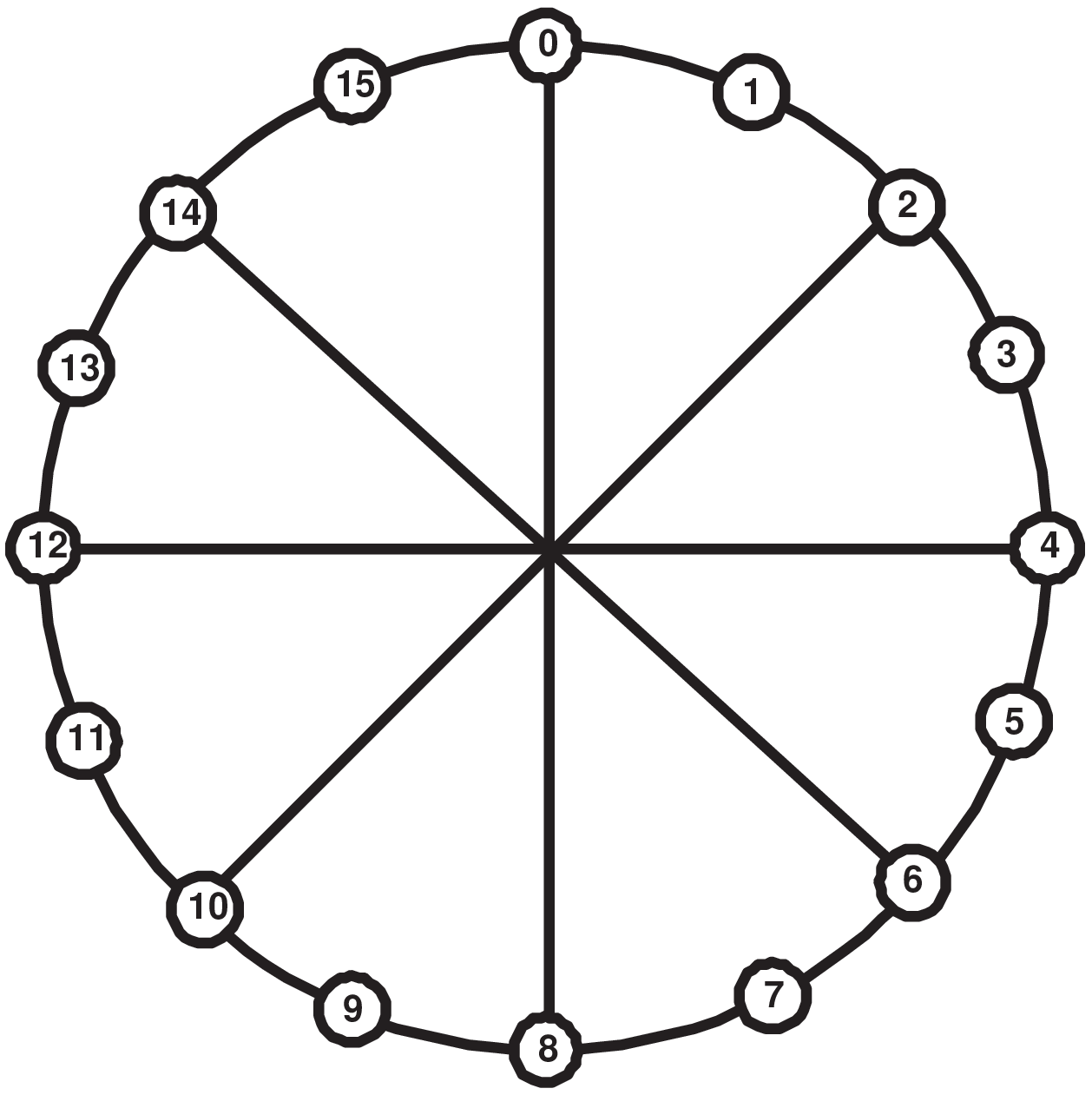}\\~\\
  Figure 3:
  The graph $G_{2}$.
 \end{center}
 Notice that
 $\chi_c(G_n)\leq 8n/(4n-1)$, since it can be verified that the
 mapping $f$ from $V(G_n)$ to $\{0,1,2,\ldots,8n-1\}$, defined below,
 is an $(8n,4n-1)$-coloring of $G_n$:
 \begin{center}
 $f(v_{2i})=2i$ and $f(v_{2i+1})=2i+1+4n($mod $8n)$ for
 $i=0,1,2,\ldots, 4n-1$.
 \end{center}
 It is also easy to check that $G_n-\{v_0,v_{4n}\}$ is bipartite,
 and deleting any vertex of $G_n$ is not enough to destroy all odd
 cycles in it. Thus $\alpha_2(G_n)=8n-2$. By Theorem
 \ref{alphat} it follows that $\chi_c(G_n)\geq
 2|V(G_n)|/\alpha_2(G_n)=
 8n/(4n-1)$. Therefore $\chi_c(G_n)=2|V(G_n)|/\alpha_2(G_n)$.

 In the following two results, we use dynamic view of circular colorings
 to extract useful information for
 $\chi_c(G)$ from
 the chromatic number of vertex's distance-2 neighborhood.
 In the sequel, we say that $\omega$ is an {\em optimal acyclic
 orientation} of $G$ with period
 $\omega_1,\omega_2,\cdots,\omega_{p_\omega}$
 if $p_\omega/m_\omega=\chi_c(G)$, and $\omega_1,\omega_2,\cdots,\omega_{p_\omega}$ is a period
 generated by $\omega$.


 \begin{theorem}
 \label{new}
 Let $H$ be a graph with vertex set
 $\{v_1,v_2,\ldots,v_n\}$.
 Let $G$ be the graph obtained from $n+1$
 disjoint graphs $H,H_1,H_2,\ldots,H_n$ by joining all vertices in $H_1,H_2,\ldots,H_n$
 to a new vertex $x$, and joining all vertices in $H_i$ to $v_i$,
 for $i=1,2,\ldots,n$. The graph $G$ is
 represented diagrammatically in Figure 4 left.
 If $H_1,H_2,\ldots,H_n$ are $t$-chromatic graphs and $\chi(H)\geq 3$,
 then $\chi_c(G)\geq t+2$.
 \end{theorem}

 \begin{proof}
 Let $\omega$ be an optimal acyclic
 orientation of $G$ with period
 $\omega_1,\omega_2,\cdots,\omega_{p_\omega}$.
 Without loss of generality, suppose that $x\in {\rm sink}(\omega_i)$.
 Consider the positive integer $s$ such that
 $$x\not \in \bigcup_{k=1}^{s}{\rm sink}(\omega_{i+k})\mbox{ and } x\in {\rm
 sink}(\omega_{i+s+1}).$$
 The fact that $x$ is fireable in $\omega_i$ and
 $\omega_{i+s+1}$ implies
 $V(H_j)\subseteq \bigcup_{k=1}^{s}{\rm sink}(\omega_{i+k})$ for $j=1,2,\ldots,n$.
 Since $H_1$ is a $t$-chromatic graph and ${\rm
 sink}(\omega_{i+1}), {\rm sink}(\omega_{i+2}), \ldots, {\rm
 sink}(\omega_{i+s})$
 are independent sets, we have $s\geq t$.
 Moreover if $s=t$ then it must be $$\mbox{$V(H_j)\bigcap {\rm sink}(\omega_{i+k})\not
 =\emptyset$,}$$
 for $j=1,2,\ldots,n$ and $k=1,2,\ldots,s$.

 Let $r$ be the chromatic number of $H$.
 We further claim that if $s=t$ then
 $$ x\not \in \bigcup_{k=2}^{t+r}{\rm sink}(\omega_{i+t+k}).$$
 To prove this, consider the positive integer $\ell$ such that
 $$x\not \in \bigcup_{k=2}^{\ell}{\rm sink}(\omega_{i+t+k})\mbox{ and } x\in {\rm
 sink}(\omega_{i+t+\ell+1}).$$
 We shall show that $\ell\geq t+r$.
 According to the previous argument,
 if $s=t$ then every vertex $y$ in $V(H_1)\cup V(H_2)\cup \ldots \cup
 V(H_n)$ becomes fireable exactly once within
 $\omega_{i+1},\omega_{i+2},\ldots,\omega_{i+t}$,
 and all vertices in $H$ are not fireable within
 $\omega_{i+1},\omega_{i+2},\ldots,\omega_{i+t}$.
 From this it follows that the digraph $\omega_{i+t+1}$ has $H_j\subseteq N^-(v_{j})$
 for $j=1,2,\ldots,n$ (see Figure 4 right).
 Since
 $$x\in {\rm sink}(\omega_{i+t+1})\bigcap {\rm
 sink}(\omega_{i+t+\ell+1}),$$
 we see that all vertices in $H_1,H_2,\ldots,H_n$
 must become fireable at least one time within
 $\omega_{i+t+2},\omega_{i+t+3},\ldots,\omega_{i+t+\ell}$.
 To ensure this, because the structure of $\omega_{i+t+1}$, we must have
 $$V(H)\subseteq \bigcup_{k=1}^\ell {\rm sink}(\omega_{i+t+k}).$$
 Suppose $\bar{\ell}\geq 2$ is the smallest integer such that
 $$V(H)\subseteq \bigcup_{k=1}^{\bar{\ell}} {\rm sink}(\omega_{i+t+k}).$$
 We see that $\bar{\ell}\geq \chi(H)=r$.
 Without loss of generality we assume $v_n\in {\rm
 sink}(\omega_{i+t+\bar{\ell}})$ and $v_n\not \in {\rm
 sink}(\omega_{i+t+k})$ for $k=1,2,\ldots,\bar{\ell}-1$.
 Since $v_n$ is an out-neighbor of every vertex of $H_n$ in $\omega_{i+t+k}$ for $k=1,2,\ldots,\bar{\ell}$
 and since
 $V(H_n)\subseteq \cup_{k=2}^{\ell} {\rm sink}(\omega_{i+t+k})$.
 We conclude that
 $$V(H_n)\subseteq \bigcup_{k=\bar{\ell}+1}^{\ell} {\rm sink}(\omega_{i+t+k}), $$
 and therefore $\ell-\bar{\ell}\geq \chi(H_n)=t$. Thus $\ell\geq t+r$.
 So we have shown that
 if $x\in {\rm sink}(\omega_i)$, $x\not \in \bigcup_{k=1}^{s}{\rm sink}(\omega_{i+k})$ and $x\in {\rm
 sink}(\omega_{i+s+1})$
 then it must be either $s\geq t+1$, or
 $$s=t\mbox{ and }x\not\in \bigcup_{k=2}^{t+r}{\rm
 sink}(\omega_{i+t+k}).$$

 Note that there are exactly $m_\omega$ integers
 $1\leq i_1< i_2< \ldots < i_{m_\omega} \leq p_\omega$ such
 that ${\rm sink}(\omega_{i_k})$ contains $x$ for
 $k=1,2,\ldots,m_\omega$.
 Let $\ell_k=i_{k+1}-i_k$ for $k=1,2,\ldots,m_\omega-1$,
 and let $\ell_{m_\omega}=p_\omega-(i_{m_\omega}-i_1)$.
 Clearly, we have $\ell_k\geq t+1$ for $k=1,2,\ldots,m_\omega$.
 Furthermore, if $\ell_k=t+1$
 then $\ell_{k+1}$ must be greater or equal than $t+r$
 (subindex addition is done modulo $m_\omega$).
 Since $r\geq 3$, it follows that
 $$p_\omega=\sum_{k=1}^{m_\omega}\ell_k\geq
 \sum_{k=1}^{m_\omega}(t+2)=m_\omega (t+2),$$
 and therefore $\chi_c(G)=p_\omega/m_\omega\geq t+2$.
 \end{proof}


 \begin{center}
 \qquad\includegraphics[scale=0.7]{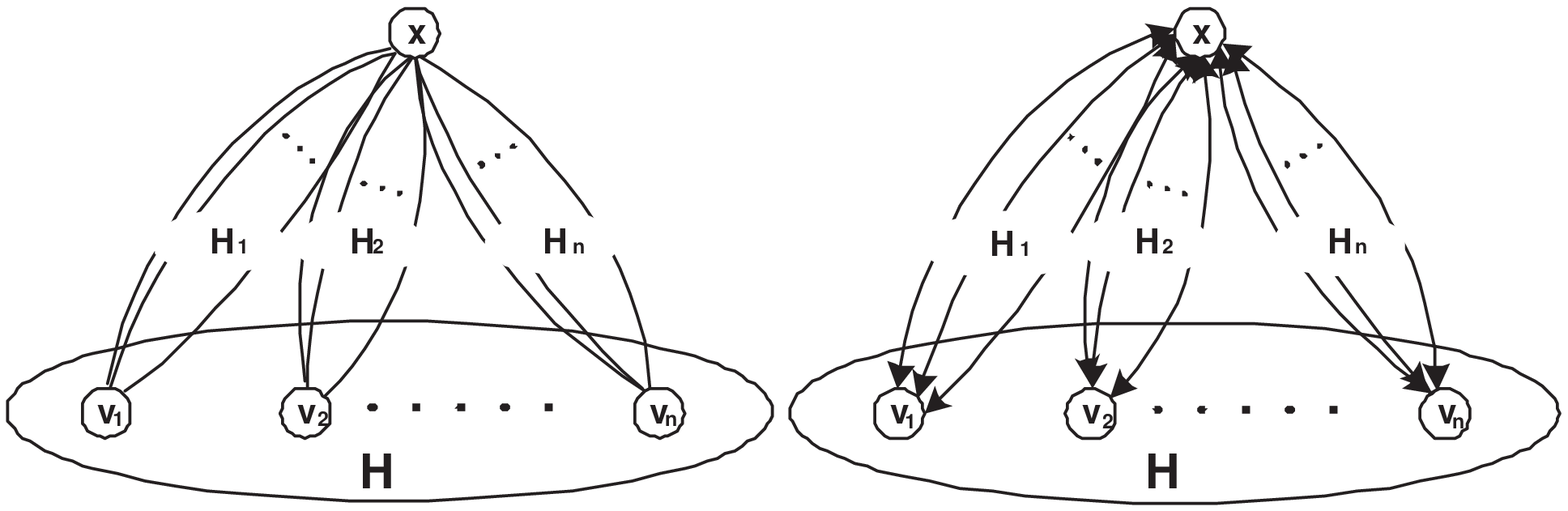}\\~\\
  Figure 4: The graph $G$ (left) and the digraph $\omega_{i+t+1}$ (right).
 \end{center}






 \begin{corollary}
 \label{d2}
 If there is a vertex $u$ in graph $G$ with $\chi(N_2(u))\geq 3$
 then $\chi_c(G)\geq 3$.
 \end{corollary}

 \begin{proof}
 Because $\chi(N_2(u))\geq 3$, graph $G$ has an odd cycle
 $$C_{2n+1}=(v_1,v_2,\ldots, v_{2n+1},v_1)$$
 such that $d_G(u,v_k)=2$ for each
 $k=1,2,\ldots,2n+1$.
 Let $ux_kv_k$ be a path in $G$ for $k=1,2,\ldots,2n+1$.
 Let $W_G$ be the subgraph of $G$ induced by
 $\cup_{k=1}^{2n+1}\{u,x_k,v_k\}$.
 Suppose graph $W$ is obtained from the odd cycle $C_{2n+1}$
 and the vertex $u$
 by adding a new vertex $h_k$ and joining $h_k$ to vertices $u$ and
 $v_k$ for all $k=1,2,\ldots,2n+1$ (see Figure 5).
 Clearly, a $(k,d)$-coloring of $W_G$ is also a $(k,d)$-coloring of
 $W$. Therefore, by Theorem \ref{new}, we have
 $$ \chi_c(G)\geq \chi_c(W_G)\geq \chi_c(W)\geq 3.$$
 \end{proof}


 \begin{center}
 \qquad\includegraphics[scale=0.7]{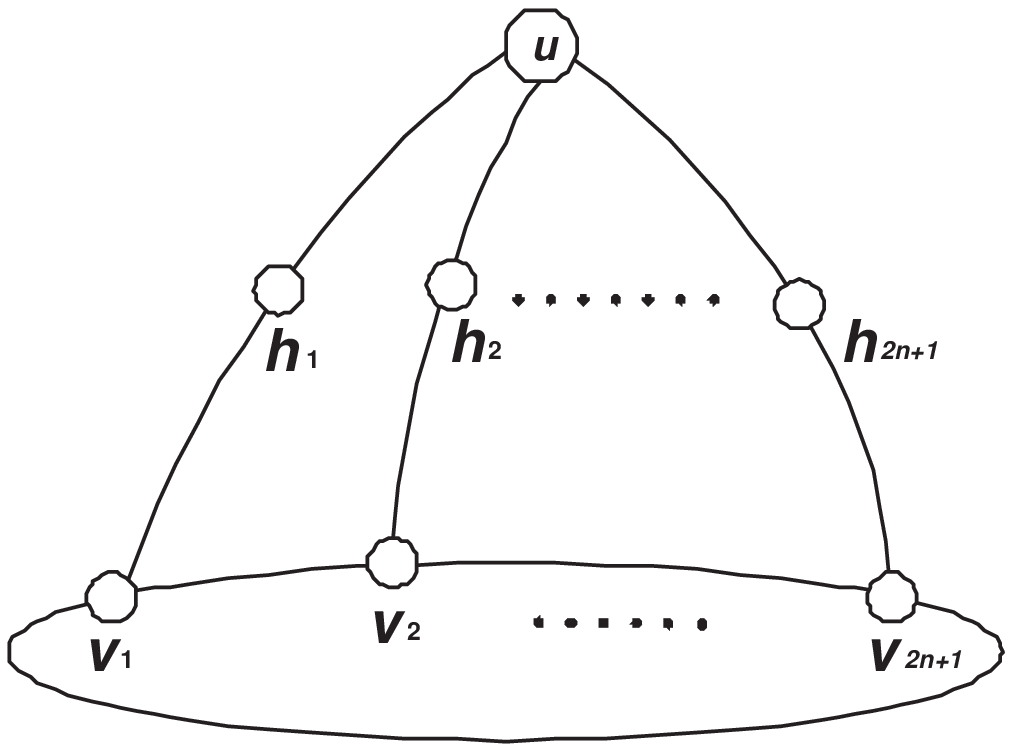}\\~\\
  Figure 5: The  graph $W$.
 \end{center}
 By using Theorem \ref{new}, we can show that graphs depicted in
 Figure 6 has $\chi_c(G_1)=3$ and $\chi_c(G_2)=4$.


 \begin{center}
 \qquad\includegraphics[scale=0.8]{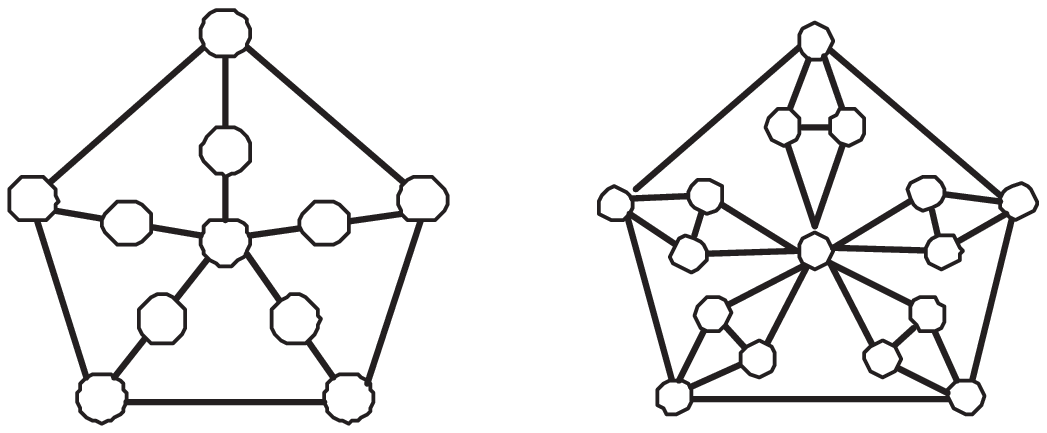}\\~\\
  Figure 6:
  Graphs $G_1$ (left) and $G_2$ (right).
 \end{center}

 By studying the behavior of the period generated by a graph's optimal
 acyclic orientation, in the following two theorems, we will establish lower bounds
 for circular chromatic number of the
 form
 $$ \chi_c(G)\geq { |V(G)|\over \alpha_1(G)-\epsilon}.  $$


 \begin{theorem}
 \label{alpha2}
 If a graph $G$ has the following four properties P1, P2,
 P3 and P4,
 then $$ \chi_c(G)\geq { |V(G)|\over \alpha_1(G)-(2/3)}. $$

 P1: $\chi(N_1(u,v))\geq 2$ for any two nonadjacent vertices $u$
 and $v$.

 P2: $|V(G)|\leq 3\alpha_1(G)-3$.

 P3: $\alpha_2(G)<2\alpha_1(G)$.

 P4: $\chi(G)\geq 3$.
 \end{theorem}

 \begin{proof}
 Let $\omega$ be an optimal acyclic
 orientation of $G$ with period
 $\omega_1,\omega_2,\cdots,\omega_{p_\omega}$.
 It suffices to consider the case where there is an index $i$ with
 $|{\rm sink}(\omega_{i})|=\alpha_1(G)$, for otherwise we have
 $$ |V(G)|m_\omega=\sum_{i=1}^{p_\omega}|{\rm sink}(\omega_{i})|
 \leq p_\omega (\alpha_1(G) -1) $$
 and hence $p_\omega/m_\omega \geq |V(G)|/(\alpha_1(G)-1)$.
 Without loss of generality, we may assume that
 $|{\rm sink}(\omega_{1})|=\alpha_1(G)$.

 Since $\alpha_2(G)<2\alpha_1(G)$,
 we see that $\langle \alpha_1(G),\alpha_1(G)\rangle$ is not a
 sub-pattern of $\omega$.
 Furthermore, $\langle \alpha_1(G),\alpha_1(G)-1,\alpha_1(G)\rangle$ is also not a
 sub-pattern of $\omega$, for otherwise the property P2 will imply
 that, for some index $i$, $$|{\rm sink}(\omega_{i})\bigcap{\rm sink}(\omega_{i+2})|\geq
 2.$$  But that is impossible since the property
 P1. Thus we have already shown that, for any index $i$,
 $$ |{\rm sink}(\omega_{i})|+
                |{\rm sink}(\omega_{i+1})|+
                |{\rm sink}(\omega_{i+2})|\leq 3\alpha_1(G)-2. $$
 Notice that property P4 guarantees that $p_\omega \geq 3$.
 By the above observations, we also have either
 $$\mbox{$|{\rm sink}(\omega_{p_\omega-1})|\leq \alpha_1(G)-1$ and
 $|{\rm sink}(\omega_{p_\omega})|\leq \alpha_1(G)-1$,}$$
 or
 $$\mbox{$|{\rm sink}(\omega_{p_\omega-1})|=\alpha_1(G)$ and
 $|{\rm sink}(\omega_{p_\omega})|\leq \alpha_1(G)-2$.}$$
 Let
 $r=p_\omega - 3 \lfloor p_\omega/3 \rfloor$.
 We have
 \begin{eqnarray*}
 |V(G)|m_\omega & = & \sum_{i=1}^{p_\omega}|{\rm
 sink}(\omega_{i})|\\
                & = & \sum_{t=1}^{\lfloor p_\omega/3 \rfloor}
                \left(
                |{\rm sink}(\omega_{3t-2})|+
                |{\rm sink}(\omega_{3t-1})|+
                |{\rm sink}(\omega_{3t})|
                \right)\\
                &      &+
                \sum_{i=p_\omega-r+1}^{p_\omega}|{\rm
                sink}(\omega_{i})|\\
                & \leq & \sum_{t=1}^{\lfloor p_\omega/3 \rfloor}
                \left( 3\alpha_1(G)-2 \right)
                +
                \sum_{i=p_\omega-r+1}^{p_\omega}
                \left( \alpha_1(G)-1 \right)\\
                & = & p_\omega (\alpha_1(G)-1)+
                {\lfloor {p_\omega/ 3} \rfloor}\\
                & \leq & p_\omega \left(\alpha_1(G)-{2/ 3} \right)
 \end{eqnarray*}
 Now we arrive at the desired inequality
 $ p_\omega/m_\omega\geq { |V(G)|\over
 \alpha_1(G)-(2/3)}$.
 \end{proof}



 \begin{theorem}
 \label{alpha1-1}
 Let $t$ be a positive integer.
 If a graph $G$ has the following five properties P1, P2,
 P3, P4 and P5,
 then $$ \chi_c(G)\geq { |V(G)|\over \alpha_1(G)-{t-1\over t}}. $$

 P1: $\chi(N_1(v))\geq t-2$ for any vertex $v$ of $G$.

 P2: $\chi(N_1(I))\geq t-1$ for any maximum independent set $I$ of
 $G$.

 P3: $|V(G)|\leq t\alpha_1(G)-t$.

 P4: Any two different maximum independent sets of $G$ intersect in
 exactly one

  vertex.

 P5: $\chi(G)\geq t$.
 \end{theorem}

 \begin{proof}
 Let $\omega$ be an optimal acyclic
 orientation of $G$ with period
 $\omega_1,\omega_2,\cdots,\omega_{p_\omega}$.

 Assume that there are two positive integers $i$, $\ell$ ($\ell\leq
 t-1$) such that
 $$|{\rm sink}(\omega_{i})|=|{\rm sink}(\omega_{i+\ell})|=\alpha_1(G),$$
 $$\mbox{$|{\rm sink}(\omega_{k})|=\alpha_1(G)-1$ for any $k\in
 [i+1,i+\ell-1]$}$$
 and
 $$\mbox{$|{\rm sink}(\omega_{i})\bigcap{\rm
 sink}(\omega_{i+\ell})|=\{u\}$
 for some vertex $u$}.$$
 Since $\chi(N_1(u))\geq t-2$, we have $\ell\geq t-1$ and hence
 $\ell=t-1$. Furthermore, since
 $$ \sum_{s=0}^{t-1}|{\rm sink}(\omega_{i+s})|=t(\alpha_1(G)-1)+2\geq |V(G)|+2, $$
 the property P4 implies that
 there are two different integers $i_1$ and $i_2$
 with $|i_1- i_2|\leq t-2$ such that
 ${\rm sink}(\omega_{i_1})\cap{\rm sink}(\omega_{i_2})\not= \emptyset$.
 However, this is impossible since $\chi(N_1(v))\geq t-2$ for any
 vertex $v$ in $G$. Thus our assumption is not true.

 Since the property P2, we note that if
 $|{\rm sink}(\omega_{i})\bigcap{\rm sink}(\omega_{j})|=\alpha_1(G)$
 for two indices $i<j$ then we must have $j\geq i+t$.
 From what we have already shown, it is clear that
 there are only two scenarios for the period
 $\omega_1,\omega_2,\cdots,\omega_{p_\omega}$:\\
 $\bullet$ Scenario 1: There is at most one index $i$ in $\{1,2,\cdots,p_\omega\}$
 with $|{\rm sink}(\omega_{i})|=\alpha_1(G)$\\
 $\bullet$ Scenario 2: If
 $|{\rm sink}(\omega_{i})|=|{\rm sink}(\omega_{j})|=\alpha_1(G)$
 for two indices $i<j$
 then either $ j\geq i+t$,
 or $ \mbox{$j=i+t-1$ with $|{\rm sink}(\omega_{k})|< \alpha_1(G)-1$ for some
 $k\in [i+1,j-1]$.} $

 In fact we have already shown that, for any index $i$,
 $$ \sum_{s=0}^{t-1}|{\rm sink}(\omega_{i+s})|\leq t(\alpha_1(G)-1)+1. $$
 As in the proof of Theorem \ref{alpha2}, we only need to consider
 the case where there is an index $i$ with
 $|{\rm sink}(\omega_{i})|=\alpha_1(G)$, and
 we may assume that
 $|{\rm sink}(\omega_{1})|=\alpha_1(G)$.
 Let $r=p_\omega - t \lfloor p_\omega/t \rfloor$.
 By the property P5 and the above observations, we have
 $$ \sum_{i=p_\omega -r+1}^{p_\omega}|{\rm sink}(\omega_{i})|
 \leq r(\alpha_1(G)-1). $$
 With these results we have
 \begin{eqnarray*}
 |V(G)|m_\omega & = & \sum_{i=1}^{p_\omega}|{\rm
 sink}(\omega_{i})|\\
                & = & \sum_{j=1}^{\lfloor p_\omega/t \rfloor}
                \left(
                \sum_{k=0}^{t-1}
                |{\rm sink}(\omega_{jt-k})|
                \right)
                +
                \sum_{i=p_\omega-r+1}^{p_\omega}|{\rm
                sink}(\omega_{i})|\\
                & \leq & \sum_{j=1}^{\lfloor p_\omega/t \rfloor}
                \left( t\alpha_1(G)-t+1 \right)
                +
                r
                \left( \alpha_1(G)-1 \right)\\
                & = & p_\omega (\alpha_1(G)-1)+
                {\lfloor {p_\omega/ t} \rfloor}\\
                & \leq & p_\omega \left(\alpha_1(G)-{(t-1)/ t}
                \right).
 \end{eqnarray*}
 Clearly we have $$ {p_\omega\over m_\omega}\geq { |V(G)|\over
 \alpha_1(G)-{t-1\over t}}.$$
 \end{proof}

 Let $H$ and $Q$ be the graphs depicted in Figure 7,
 where $Q$ is obtained from the Petersen graph
 by deleting one vertex.



 \begin{center}
 \qquad\includegraphics[scale=0.7]{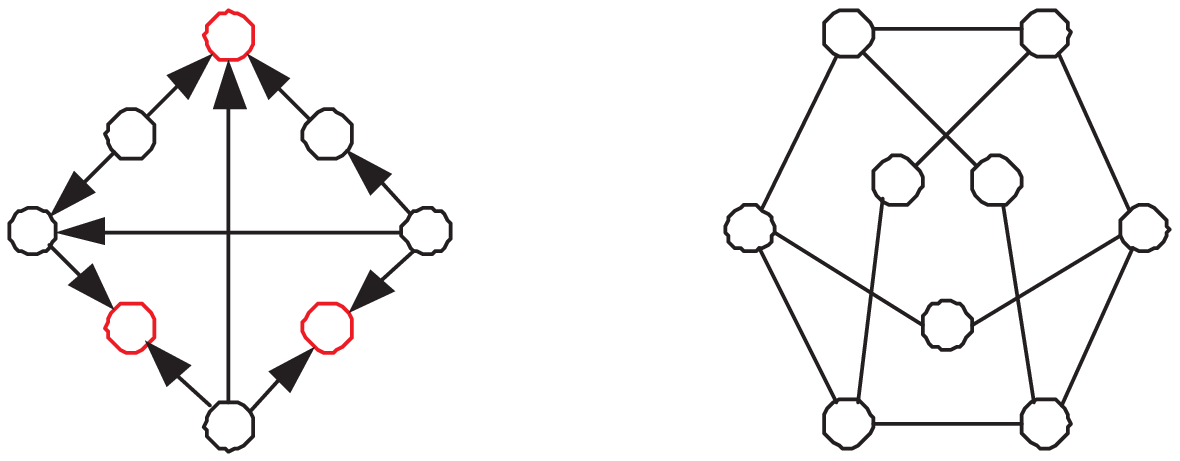}\\~\\
  Figure 7:
  Digraph $(H,\omega)$ (left) and graph $Q$ (right).
 \end{center}



 \begin{example}
 \label{H}
 $\chi_c(H)=8/3$.
 \end{example}

 \begin{proof}
 It can be checked easily that
 the acyclic orientation $\omega$ of $H$ (as depicted in Figure 7)
 has $p_\omega /m_\omega =8/3$,
 and hence $\chi_c(H)\leq 8/3$ by Corollary \ref{pm}.
 Moreover, by Theorem \ref{alphat},
 one has $\chi_c(H)\geq 2|V(H)|/\alpha_2(H)= 16/\alpha_2(H)\geq 8/3$.
 \end{proof}



 \begin{example}
 \label{Q}
 $\chi_c(Q)=3$.
 \end{example}

 \begin{proof}
 From
$\chi_c(Q)\in \{ {k \over d} : \mbox{$k\leq |V(Q)|$, $d\leq
\alpha(Q)$ and }
   { |V(Q)| \over \alpha(Q) } \leq {k \over d} \leq \chi(Q) \}$,
 we have $\chi_c(Q)\in\{  {5\over 2},{8\over 3},3  \}$.
 Since $\alpha_1(Q)=4$, $\alpha_2(Q)<8$ and $\chi(Q)=3$,
 we can easily check that
 the graph $Q$ satisfies all the
 properties stated in the Theorem \ref{alpha2}.
 It follows that $\chi_c(Q)\geq {|V(Q)|\over
 \alpha_1(Q)-(2/3)}=27/10>8/3$, which completes the proof.
 \end{proof}



 \begin{example}
 \label{PL}
 Suppose $P_L$ is the line graph of the Petersen graph.
 Then
 $\chi_c(P_L)=11/3$.
 \end{example}

 \begin{proof}
 The acyclic orientation $\omega$ of $P_L$ (depicted in Figure 8) has
 $p_{\omega} / m_{\omega}=11/3$,
 and hence $\chi_c(P_L)\leq 11/3$.
 Since $\alpha_1(P_L)=5$ and $\chi(P_L)=4$, similar to the proof of Example \ref{Q},
 we have, by \cite{survey},
 $\chi_c(P_L) \in
 \{ 3,{13\over 4}, {10\over 3}, {7\over 2}, {11\over 3} \}$.
 We note that the subgraph left by deleting the edges of a perfect matching
 from
 the Petersen graph contains an odd cycle.  We also see that any
 two different maximum matchings of the Petersen graph intersect in
 exactly one edge. Therefore $P_L$ satisfies all the
 properties stated in Theorem \ref{alpha1-1} for $t=4$.
 So we have
 $$ \chi_c(P_L)\geq { |V(P_L)|\over \alpha_1(P_L)-{3\over 4}}=
 {15\over 5-{3\over 4}}={60\over 17}>{7\over 2}. $$
 Thus it must be $\chi_c(P_L)= 11/3$.
 \end{proof}



 \begin{center}
 \qquad\includegraphics[scale=0.75]{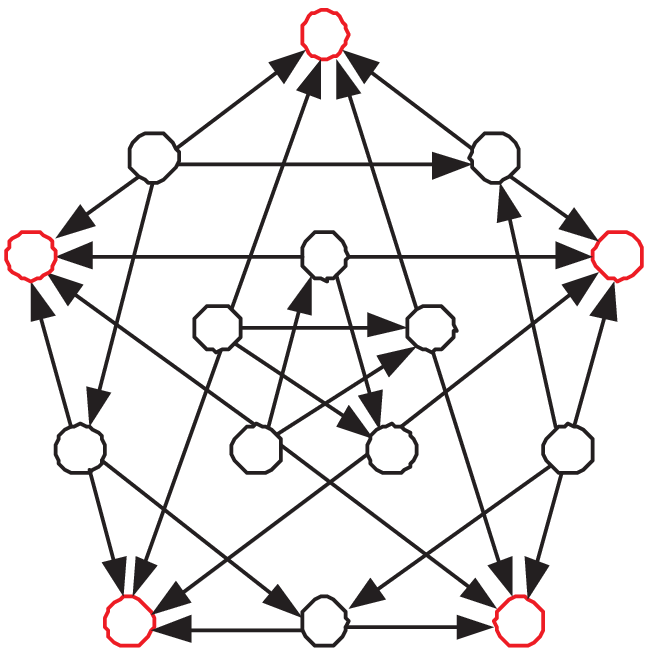}\\~\\
  Figure 8: Digraph $(P_L, \omega)$.
 \end{center}

%
%

\end{document}